\documentclass[12pt]{amsart} 
\usepackage{verbatim} \usepackage{varioref} \usepackage{array}
\usepackage{amssymb}
\usepackage[mathscr]{eucal}
\usepackage{amsthm} \usepackage{topmacs} \usepackage[all]{xy}

\CompileMatrices

\SelectTips{cm}{}
\begin{document}

\title[m-structures]{m-structures determine integral homotopy type}
\author{Justin R. Smith}
\date{June 9, 1998}

\begin{abstract}
This paper proves that the functor $\mathscr{C}(*)$ that sends
pointed, simply-connected CW-complexes to their chain-complexes
equipped with diagonals and iterated higher diagonals, determines
their integral homotopy type --- even inducing an equivalence of
categories between the category of CW-complexes up to homotopy
equivalence and a certain category of chain-complexes equipped with
higher diagonals.  Consequently, $\mathscr{C}(*)$ is an algebraic
model for integral homotopy types similar to Quillen's model of
rational homotopy types.  For finite CW complexes, our model is
finitely generated. 

Our result implies that the geometrically induced diagonal map with
all ``higher diagonal'' maps (like those used to define Steenrod
operations) collectively determine integral homotopy type.
\end{abstract}
\maketitle

\section{Introduction}

This paper forms a sequel to \cite{Smith:1994}. That paper developed 
the theory of m-coalgebras and defined a functor $\mathscr{C}(*)$ 
that associated canonical m-coalgebras to  semi-simplicial 
complexes. 

\newdir{ >}{{}*!/-5pt/@{>}}
\def\contractright#1{\ar@<.5ex>@{ >-}[r]^-{#1}\ar@<.5ex>@{~>>}[r];[]^-{#1'}}
\def\contractleft#1{\ar@<-.5ex>@{ >-}[l]_-{#1}\ar@<-.5ex>@{~>>}[l];[]_-{#1'}}
\def\contractup#1{\ar@<.5ex>@{ >-}[u]^-{#1}\ar@<.5ex>@{~>>}[u];[]^-{#1'}}

Our main result is:
\medskip

{\bf Corollary~\zref{cor:mainresult}:} {\it
The functor  (defined in 4.2 on page 30 of \cite{Smith:1994})
\[
\mathscr{C}(*):\underline{\mathrm{Homotop}}_0\to \hat{\mathfrak{M}}
\]
(see \zref{def:bigmcat} for the definition of $\hat{\mathfrak{M}}$)
defines an equivalence of categories, where
$\underline{\mathrm{Homotop}}_0$ is the category of pointed,
simply-connected CW-complexes and continuous maps, in which homotopy
equivalences have been inverted (i.e., it is the category of fractions
by homotopy equivalences).}
\medskip

This, of course, implies the claim made in the title --- that
m-structures determine integral homotopy type.

From the beginning, it has been a central goal of homotopy theory (and
algebraic topology in general) to develop tractable models for spaces
and mappings.  The early models were combinatorial, including
simplicial or semi-simplicial complexes, chain-complexes, DGA algebras
and coalgebras and so on.  These models tended to fall into two
classes:
\begin{itemize}

\item powerful, but computationally intractable (i.e., minimal models,
chain-complexes of free rather than abelian groups, etc.)

\item weak (chain-complex) but well-behaved.
\end{itemize}

The first major breakthrough came with the work of Quillen in 
\cite{Quillen:1969}, in which he simplified the problem by focusing 
on rational homotopy types. Rationalizing eliminates much of homotopy 
theory's complexity by killing off cohomology operations, like 
Steenrod operations.  Quillen was able to create a complete and 
faithful model of rational homotopy theory --- co-commutative 
DGA-coalgebras over $\rationals$. 

This paper is the outcome of a research program of several years 
duration. One of the main goals of this program was to understand
the coproduct or cup-product structure of the total space of a 
fibration.  In order to accomplish this, it was necessary to compute 
a topological coproduct on the cobar construction and on the 
canonical acyclic twisted tensor with fiber a cobar construction.

Although the cobar construction is defined for DGA-coalgebras, 
computing a ``geometric'' coproduct on the cobar construction requires 
more than the mere coproduct.  I quickly realized that various 
cohomology operations entered into the cobar construction's coproduct.  
It was necessary to equip the chain complex of a space with diagonals 
and higher diagonals defined on the {\it chain level\/} (rather than 
on cohomology with coefficients in a finite field).

These higher coproducts satisfy a complex web of relationships I call 
{\it coherence conditions.\/} In \cite{Smith:1994}, I developed an 
algebraic device called an  m-coalgebra over a formal coalgebra to 
encapsulate these relations.  A referee of \cite{Smith:1994} pointed 
out that formal coalgebras had been defined and studied before under 
the name {\it operad.} 

This research had a gratifying outcome: A coherent 
m-coalgebra's cobar construction not only has a computable coproduct; 
it comes equipped with a well-defined and geometrically valid {\it 
m-coalgebra structure\/} (although the coherence condition must be 
weakened slightly).  It, consequently, becomes possible to {\it 
iterate\/} the cobar construction.  A side-effect was an explicit 
procedure for computing geometric m-coalgebra structures on the total 
space of a fibration (represented by a twisted tensor product).

This suggested to me a {\it possibility\/} of characterizing integral 
homotopy theory: If one can compute coproducts (and higher
coproducts) on fibrations, one can in principal compute fibrations 
over fibrations, and so on. This suggested the possibility of 
purely algebraic computations of {\it Postnikov towers} --- possibly 
along the lines of Sullivan in \cite{Sullivan.mm:1977}.

The present paper is the result.

In 1985, Smirnov proved a result similar to ours in
\cite{Smirnov:1985} --- showing that a functor whose value is a
certain comodule over a certain operad determines the integral
homotopy type of a space. The operad and comodule in question were
uncountably generated in all dimensions and in the simplest case.

In contrast, our functor is finitely generated in all dimensions for
finite simplicial complexes.  Although it is considerably more complex
than the co-commutative coalgebras Quillen derived, it is highly
unlikely one can get away with something much simpler: {\it all\/} of
our functor appears nontrivially in even the {\it coproduct\/} of a
cobar construction.

At this point, I feel it is appropriate to compare and contrast my 
results with work of Michael Mandell.  In \cite{mandell:1998-1}, he 
proved

\medskip
\begingroup
{\bf Main Theorem.}\it
The singular cochain functor with coefficients in 
$\bar{\integers}_{p}$ induces a contravariant equivalence from the 
homotopy category of connected nilpotent $p$-complete spaces of 
finite $p$-type to a full subcategory of the homotopy category of 
$E_{\infty}$ $\bar{\integers}_{p}$-algebras.
\endgroup
\medskip

Here, $p$ denotes a prime and $\bar{\integers}_{p}$ the algebraic 
closure of the finite field of $p$ elements.  $E_{\infty}$-algebras 
are defined in \cite{Kriz-May} --- they are modules over a suitable 
operad.

At first glance, it would appear that his results are a kind of dual 
to mine: He characterizes nilpotent $p$-complete spaces in terms of 
$E_{\infty}$ $\bar{\integers}_{p}$-algebras.  This is not the case, 
however.  A complete characterization of nilpotent $p$-complete spaces 
does not lead to one of integral homotopy types:  One must somehow 
know that $p$-local homotopy equivalences patch together.
Consequently, his results do not imply mine.

The converse statement is also true: My results do not imply his.

My results in \cite{Smith:1994} imply that all
the primes ``mix'' when one studies algebraic properties of homotopy
theory (for instance the $p$-local structure of the cobar construction
of a space depend on the $q$-local structure of the space for all
primes $q\ge p$).  This is intuitively clear when considers the
composite $(1\otimes \Delta)\circ \Delta$ (iterated coproducts) and
notes that $\integers_{2}$ acting on both copies of $\Delta$ give rise
to elements of the symmetric group on 3 elements.

Consequently, a characterization of integral homotopy does not lead to
a $p$-local homotopy theory: In killing off all primes other than $p$,
one also kills off crucial information needed to compute the cobar
construction of a space.  

In \cite{mandell:1998-1}, Dr.~Mandell proved that one {\it must\/}
pass to the algebraic closure of $\integers_p$ to get a
characterication of $p$-complete homtopy theory.  I conjecture that,
in passing to the algebraic closure, one kills off additional data
within the homotopy type --- namely the data that depends on larger
primes. Consequently, one restores algebraic consistence to the
theory, regaining the ability to characterize local homotopy types.

I am indebted to Jim Stasheff for his encouragement and to Michael 
Mandell for pointing out errors and inconsistencies in an earlier 
version of this paper.

\section{Definitions and preliminaries}

We recall a few relevant facts from \cite{Smith:1994}.

\begin{definition} \label{d:koszul} If $f:C_1\to D_1$, $g:C_2\to D_2$ are
maps, and $a\otimes b \in C_1\otimes C_2$ (where $a$ is a homogeneous
element), then $(f\otimes g)(a\otimes b)$ is defined to be
$(-1)^{\deg(g)\cdot \deg(a)}f(a)\otimes g(b)$. \end{definition} 

\begin{remarks} \rem This convention simplifies many of the common
expressions that occur in homological algebra --- in particular it
eliminates complicated signs that occur in these expressions. For
instance the differential, $\partial_\otimes$, of the tensor product
$C\otimes D$ is just $\partial_C\otimes 1+1\otimes \partial_D$.

\rem  Throughout this entire paper we will follow the convention
that group-elements act on the left. Multiplication of elements of
symmetric groups will be carried out accordingly --- i.e.
$\begin{pmatrix} 1&2&3&4\\ 2&3&1&4 \end{pmatrix} *\begin{pmatrix}
1&2&3&4\\ 4&3&2&1 \end{pmatrix} =$ result of applying $\begin{pmatrix}
1&2&3&4\\ 2&3&1&4 \end{pmatrix}$ first and then $\begin{pmatrix}
1&2&3&4\\ 4&3&2&1 \end{pmatrix}$. The product is thus $\begin{pmatrix}
1&2&3&4\\ 4&3&1&2 \end{pmatrix}$.

\rem Let $f_i$, $g_i$ be maps. It isn't hard to verify that the
Koszul convention implies that $(f_1\otimes g_1)\circ (f_2\otimes g_2) =
(-1)^{\deg(f_2)\cdot\deg(g_1)}(f_1\circ f_2\otimes g_1\circ g_2)$.

\rem  We will also follow the convention that, if $f$ is a map
between chain-complexes, $\partial f = \partial\circ f -
(-1)^{\deg(f)}f\circ \partial$. The compositions of a map with boundary
operations will be denoted by $\partial\circ f$ and $f\circ \partial$ ---
see \cite{Gugenheim:1960}. This convention clearly implies that
$\partial(f\circ g) = (\partial f)\circ g + (-1)^{\deg(f)}f\circ
(\partial g)$. We will call any map $f$ with $\partial f=0$ a chain-map.
We will also follow the convention that if $C$ is a chain-complex and
${\susp}:C \to \Sigma C$ and $\desusp:C \to \Sigma^{-1}C$ are,
respectively, the suspension and desuspension maps, then ${\susp}$ and
$\desusp$ are both chain-maps. This implies that the boundary of $\Sigma C$
is $-{\susp}\circ \partial C\circ \desusp$ and the boundary of
$\Sigma^{-1}C$ is $-\desusp\circ \partial C\circ {\susp}$.

\rem \label{r:koszul.5}  We will use the
symbol $T$ to denote the transposition operator for tensor products of
chain-complexes $T:C\otimes D \to D\otimes C$, where $T(c\otimes d) =
(-1)^{\dim(c)\cdot \dim(d)}d\otimes c$, and $c \in C$, $d \in D$.
\end{remarks}

\begin{definition} \label{d:fcoalg} Let $\{U_n\}$ denote a sequence of
differential graded $\integers$-chain-complexes with preferred
$\integers$-bases, $\{b_\alpha\}$, with $n$ running from $1$ to $\infty$.
This sequence will be said to constitute an {\it operad} with
$U_n$ being the component of rank $n$ if:
\begin{quotation} \noindent
given $\integers$-basis elements, $S_1$ and $S_2$, the following
(possibly distinct) composites are defined: $\{S_1 \circ_k S_2\}$, where
$1 \le k \le \rank(S_2)$ and all are defined to have rank equal to
$\rank(S_1) + \rank(S_2) - 1$ and degree equal to $\dim(S_1) +
\dim(S_2)$. These composition operators are subject to the following
identities:
\begin{enumerate}
\item $(S_1 \circ_i S_2) \circ_j S_3 = S_1
\circ_{i+j-1} S_2 \circ_j S_3$;
\item if $j<i$ then $S_1
\circ_{i+\rank(S_2)-1} S_2\circ_j S_3 = S_2\circ_j S_1 \circ_i S_3$
\end{enumerate} \end{quotation} 

The differential $\partial:U \to U$: \begin{enumerate} \item preserves
rank; \item imposes the following additional condition on composition
operations $\partial(S_1 \circ_i S_2) = \partial S_1 \circ_i S_2 +
(-1)^{\dim(S_1)}S_1 \circ_i \partial S_2$. \end{enumerate}
\end{definition} 

\begin{remarks} \rem Multiple compositions are assumed to be
right-associative unless otherwise stated --- i.e. $S_1 \circ_i S_2
\circ_j S_3 = S_1 \circ_i (S_2 \circ_j S_3)$.

\rem  An operad will be called {\it unitary} if it
contains an identity element with respect to the composition-operations
$\{ \circ_i\}$. This will clearly have to be an element of rank 1 and
degree 0. 

\rem Our definition of an operad in the category of DGA algebras is 
slightly different from the standard one given in  \cite{Kriz-May}.
It is a simple exercise to see that the two definitions are 
equivalent:
The fundamental degree-$n$ operation of an operad, $Z$, (in the 
standard definition)
is a $n+1$-linear map
\[
Z_{i_{1}}\otimes Z_{i_{2}}\otimes\cdots\otimes Z_{i_{n}}\otimes 
Z_{n}\to Z_{i_{1}+\cdots +i_{n}}
\]
which is simply an $n$-fold iteration of our ``higher'' compositions:
\[
z_{1} \circ_{1} z_{2}\circ_{2}\cdots z_{n}\circ_{n} b
\]
where $b\in  Z_{n}$ and $z_{j}\in Z_{i_{j}}$.

Our notation lends itself to the kinds of computations we want to do.
\end{remarks}

\begin{definition} Let $A$ and $B$ be operads. A morphism $f: A
\to B$ is a morphism of the underlying chain-complexes, that preserves
the composition operations. \end{definition}

Now we give a few examples of operads:

\begin{definition} The {\it trivial} operad, denoted $I$, is
defined to have one basis element $\{b_i\}$ for all integers $i \ge 0$.
Here the rank of $b_i$ is $i$ and the degree is 0 and the these elements
satisfy the composition-law: $b_i\circ_\alpha b_j = b_{i+j-1}$,
regardless of the value of $\alpha$, which can run from $1$ to $j$. The
differential of this formal coalgebra is identically zero.
\end{definition} 

\begin{remark} This is clearly a unitary operad --- the
identity element is $b_1$. \end{remark}

\begin{definition} Let $C_1$ and $C_2$ be operads. Then
$C_1\otimes C_2$ is defined to have: \begin{enumerate} 
	
\item component of $\rank i = (C_1)_i\otimes (C_2)_i$, where $(C_1)_i$
and $(C_2)_i$ are, respectively, the components of $\rank i$ of of $C_1$
and $C_2$;
	
\item composition operations defined via $(a\otimes b) \circ_i (c\otimes
d) = (-1)^{\dim(b)\dim(c)}(a \circ_i c\otimes b \circ_i d)$, for $a, c
\in C_1, b, d \in C_2$.
\end{enumerate} \end{definition}

\begin{definition} Let $C$ be a DGA-module with augmentation $\epsilon:C
\to \integers$, and with the property that $C_0=\integers$. Then the {\it
endomorphism operad} of $C$, denoted $\mathbf{P}(C)$ is defined to be the
operad with:

\begin{enumerate}
\item component of $\rank i = \homz(C,C^i)$, with the differential
induced by that of $C$ and $C^i$. The dimension of an element of
$\homz(C,C^i)$ (for some $i$) is defined to be its degree as a map.

\item  The $\integers$-summand is generated by one element, $e$,  of rank
0. \end{enumerate}

Let $s_1 \in \homz(C,C^i)$ and $s_2 \in \homz(C,C^j)$ be elements of rank
$i$ and $j$, respectively, where $i, j \ge 1$. Then the composition
$s_1 \circ_k s_2$, where $1 \le k \le j$, is defined by: $s_1 \circ_k s_2
= \underbrace{1\otimes \cdots\otimes s_1\otimes \cdots\otimes
1}_{\text{$k^{\mathrm{th}}$ position}} \circ s_2:C \to C^{i+j-1}$. The
composition $e\circ_k s_2$ is defined in a similar way, by identifying
$e$ with the augmentation map of $C$ --- it follows that $e\circ_k s_2
\in \homz(C,C^{j-1})$, as one might expect.

The canonical subcomplex $\homz(C,C^i)$ of elements of rank $i$, is
equipped with a natural $S_i$-action --- it is defined by permutation of
the factors of the target, $C^i$. \end{definition}

\begin{remarks} \rem This is a unitary operad --- its
identity element is the identity map $\mathrm{id} \in \homz(C,C)$.

\rem In general, operads model structures like the iterated coproducts 
that occur in the endomorphism operad.  We will use operads as an 
convenient algebraic framework for defining other constructs that have 
topological applications.
\end{remarks}

\begin{prop} Let $C$ be a DGA-module. Co-associative coalgebra structures
on $C$ can be identified with morphisms $f:I\to \mathbf{P}(C)$, the the
trivial operad to the endomorphism of $C$. \end{prop}

We now define a very important operad --- the {\it symmetric
construct}. It models the formal behavior of $\{\homz(C,C^n)\}$ in which
each $C^n$ is equipped with an action of $S_n$ that permutes the factors
of $C$. 

The symmetric construct will be denoted $\mathfrak{S}$. Its components are
$\{\rs n\}_{n\in\integers^+}$, where:
\begin{enumerate}
\item $S_n$
denotes the symmetric group on $n$ objects;
\item $\rs n$ denotes the
bar-resolution of $\integers$ over $\zs n$;
\end{enumerate}

Here we follow the convention that $\rs 0= \rs 1 = \integers$,
concentrated in dimension 0. Pure elements of $\mathfrak S$ are canonical
basis elements of $\rs n$ for all values of $n$, or the generator 1 of
the $\integers$-summand (by canonical basis elements, we mean elements
of the form $[g_1|\dots|g_k] \in \rs n$).

See \S~2 of \cite{Smith:1994} for a detailed description of the 
composition operations of $\mathfrak S$. 

We are now in a position to define {\it m-structures} 

\begin{definition} \label{mstructdef} Let $C$ be a chain-complex with
$H_0(C) = \integers$. Then: \begin{enumerate} 

\item An {\it m-structure} on $C$ is defined to be a sequence of chain
maps $\ncmap{C}_n:C \to \homzs{n}(\ncoord{C}_n, C^n)$, where $\mathfrak{R} =
\{\ncoord{C}_n\}$ is some f-resolution, and $n$ is an integer that
 satisfies $0\le n<\infty$. We assume that: \begin{enumerate}

\item the composite $e_1\circ f_1:C \to C^1$, is the identity map of $C$;

\item \label{item:pointed}
and the composite $e_0\circ f_0:C \to C^0 = \integers$ coincides
with the augmentation of $C$;

\item For any $c \in C$, at most a finite number of the
$\{\ncmap{C}_n(c)\}$ are nonzero. Here $C^n$ is equipped with the
$S_n$-action that permutes the factors.

\item The adjoint will be denoted
$\widetilde{\ncmap{C}}_n:\ncoord{C}_n\otimes C \to C^n$, and is defined by
$\widetilde{\ncmap{C}}_n(r\otimes c) =
(-1)^{\dim(r)\cdot\dim(c)}\ncmap{C}_n(c)(r)$, where $r \in
\ncoord{C}_n$ and $c \in C$. With this definition in mind, we require
$\widetilde{\ncmap{C}}_n(\ncoord{C}_n\otimes C(k)) \subseteq C(k)^n$,
where $C(k)$ is the $k$-skeleton of $C$. \end{enumerate}

\item \label{ite:coherentdia} An m-structure will be called {\it
weakly-coherent} if the {\it adjoint maps} fit into 
commutative diagrams:
$$
\xymatrix@d@R+20pt{
\ncoord{C}_n\otimes\ncoord{C}_m\otimes C\ar[r]^-{\circ_i}&
\ncoord{C}_{n+m-1}\otimes C\ar[r]^-{\widetilde{\ncmap{C}}_{n+m-1}}&
C^{n+m-1}\\
\ncoord{C}_n\otimes\ncoord{C}_m\otimes
C\ar[r]_-{1\otimes\widetilde{\ncmap{C}}_m}
\ar[u]^-{1_{\ncoord{C}_n\otimes\ncoord{C}_m\otimes C}}&
\ncoord{C}_n\otimes C^m\ar[r]_-{V_{i-1}}&
C^{i-1}\otimes\ncoord{C}_n\otimes C\otimes C^{m-i}
\ar[u]_-{1\otimes\dots\otimes\widetilde{\ncmap{C}}_n\otimes\dots\otimes1}
}
$$
for all $n, m
\ge 1$ and $1 \le i \le m$.
Here $V:\ncoord{C}_n\otimes C^m\to C^{i-1}\otimes\ncoord{C}_n\otimes 
C\otimes C^{m-i}$ is the map that shuffles the factor $\ncoord{C}_n$ to the 
right of $i-1$ factors of $C$.

\item An m-structure $\{\ncmap{C}_n:C \to \homzs{n}(\ncoord{C}_n,
C^n)\}$, will be called {\it strongly coherent} (or just {\it
coherent}) if it is weakly coherent, and $\ncoord{C} = \mathfrak{S}$.
\end{enumerate}

A chain-complex, $C$, equipped with an m-structure will be called an
{\it m-coalgebra}. The maps $\ncmap{C}_n:C \to \homzs{n}(\ncoord{C}_n, C^n)$,
where $n$ is an integer such that $0\le n<\infty$, will be called the
{\it structure maps} of $C$. \end{definition}

\begin{remarks} \rem \label{rem:coherentprop} If $C$ is an incoherent
m-coalgebra we may, without loss of generality, assume that $\ncoord{C}
= \mathfrak{S}$, since the contracting homotopy, $\Phi$, that is packaged
with $\ncoord{C}$, allows us to construct a unique sequence of chain-map
$\mathfrak{S}_n=\rs n \to
\ncoord{C}_n$, for $n$ an integer such that $0\le n<\infty$. We then
compose the structure maps of the original m-coalgebra with the induced
natural transformation $\homzs{n}(\ncoord{C},*) \to 
\homzs{n}(\mathfrak{S},*)$, to
get the structure maps of the modified m-coalgebra. 

\rem An m-coalgebra can be given the following interpretation: The
adjoint isomorphism allows us to regard the structure maps as a family
of $S_n$-equivariant chain-maps $\widetilde{\ncmap{C}}_n:\rs n\otimes
C \to C^n$. The map $\widetilde{\ncmap{C}}_2:\rs 2\otimes C \to C^2$,
restricted to $[~]\otimes C$, defines a kind of coproduct on $C$,
called the {\it underlying coproduct} of the m-coalgebra. Define $\D_a
= \widetilde{\ncmap{C}}_i(a\otimes *):C \to C^i$.  These maps will be
called the {\it higher-coproducts} associated with the m-coalgebra.
The map $\D_{[(1,2)]}:C \to C^2$ defines a chain-homotopy between
$\Delta = \D_{[~]}$ and $T\circ \Delta$, where $T$ is the
transposition map defined in \zref{r:koszul.5}.

\rem The basic definitions can be stated in terms of {\it operads} in
the category of graded differential modules. Operads were originally
defined in terms of topological spaces by May in \cite{May:1972} and
this concept was extended to DG-modules by Smirnov in
\cite{Smirnov:1985}.  Essentially: \begin{enumerate}

\item  the operad $\mathfrak{S}$, constitutes an operad, and 

\item a coherent m-coalgebra is a {\it comodule over} this operad, in
the sense of \S~3 of \cite{Smirnov:1985}.
\end{enumerate} 

\rem My original definition of an m-coalgebra regarded a coherent
m-structure as a morphism of operads $\mathfrak{S} \to \powconst{C}$,
and a weakly coherent m-structure as a morphism $\ncoord{C}\to
\powconst{C}$.  Although this definition has the advantage of being
much more elegant than the one given above it doesn't lend itself to
effective computation unless $C$ is finitely generated as a
$\integers$-module --- this means:

\begin{enumerate}
\item $C_i \ne 0$ for at most a finite number of values of $i$;

\item each of these nonzero $C_i$ is, itself, finitely generated as a
$\integers$-module.
\end{enumerate}

\rem The definition of weak coherence of an m-structure can be
re-stated in terms of the maps $\{\ncmap{C}_n\}$ themselves, rather
than their {\it adjoints} $\{\widetilde{\ncmap{C}}_n\}$.  An
m-structure is weakly coherent if and only if the diagram in figure
2.2.2 on page 22 of \cite{Smith:1994} commutes for all integers $n$
such that $0\le n<\infty$.  In this diagram, the map $V'_i$ represents
the composite
\begin{multline}
\homzs{n}(\ncoord{C}_n,
C^{i-1}\otimes\homzs{m}(\ncoord{C}_m,C^m)\otimes C^{n-i})\\
\xymatrix{\ar[r]^-{i_{1}}&
\homzs{n}(\ncoord{C}_n,
C^{i-1}\otimes\homzs{m}(\ncoord{C}_m,C^m)\otimes C^{n-i})}\\
\xymatrix@C+20pt{\ar[r]^-{\homzs{n}(1,i_2)}&
\homzs{n}(\ncoord{C}_n,
C^{i-1}\otimes\homzs{m}(\ncoord{C}_m,C^m)\otimes C^{n-i})}\\
\xymatrix{\ar[r]&
\homz(\ncoord{C}_n\otimes\ncoord{C}_m,
C^{n+m-1})}
\end{multline} where $i_1$ and $i_2$ are inclusion mappings of the
$\homzs{n}$-functors in the respective $\homz$-groups. We are also
including $\homzs{i}(*,*)$ in $\homz(*,*)$, by simply forgetting that the
elements are $\zs{i}$ linear.
\end{remarks}
\begingroup
\fontsize{10}{10pt}\selectfont
$$
\xymatrix@C+20pt{
C\ar[r]^{\ncmap{C}_n}
\ar[dd]_{f_{n+m-1}}&
\homzs{n}(\ncoord{C}_n,C^n)\ar[d]|-{\homz(1,1\otimes\dots\otimes f_m
\otimes\dots\otimes1)}\\
&G\ar[d]^{V'_i}\\
\homzs{n+m-1}(\ncoord{C}_{n+m-1},C^{n+m-1})\ar[r]_{\homz(\circ_i,1)}
&\homz(\ncoord{C}_n\otimes\ncoord{C}_m,C^{n+m-1})}
$$
\endgroup
where $G=\homzs{n}(\ncoord{C}_n,C^{i-1}\otimes
\homzs{m}(\ncoord{C}_m,C^m)\otimes C^{n-i})$.
This diagram means that the composition-operations in the coordinate
coalgebra correspond to actual compositions of the adjoint maps.

Coherence of an m-structure implies a number of identities involving 
compositions of higher coproducts.  For instance, $\D_{[(1,2)]}\otimes 
1\circ \D_{[(1,2)]} =\D_{[(1,2)]\circledast \Tmap{2,1}[(1,2)]} = 
\D_{[(1,3,2)]\circledast [(1,2)]} = 
\D_{[(1,3,2)|(1,2)]-[(1,2)|(1,2,3)]} = \D_{[(1,3,2)|(1,2)]} - 
\D_{[(1,2)|(1,2,3)]}$.  In fact, we can translate {\it any\/} formula 
involving compositions of higher-coproducts into one without 
compositions involving elements of the $\{\rs n\}$.

\begin{prop} \label{mcoalgprod} Let $\mathfrak{R}_1 = 
\{\mathfrak{R}_{1,n}\}$ and $\mathfrak{R}_2 = \{\mathfrak{R}_{2,n}\}$ 
be f-resolutions, and let $C_1$ and $C_2$ be chain-complexes.  Then 
there exists a natural transformation of functors 
$\mathfrak{E}_n:\homzs{n}(\mathfrak{R}_{1,n},C_1^n)\otimes 
\homzs{n}(\mathfrak{R}_{2,n},C_2^n)\to \homzs{n}(\mathfrak{R}_{1,n} 
\otimes\mathfrak{R}_{2,n},{(C_1\otimes C_2)}^n)$, for all $n$.  
\end{prop}

\begin{remark} If $u \in \homzs{n}(\mathfrak{R}_{1,n},C_1^n)$, $v \in 
\homzs{n}(\mathfrak{R}_{2,n},C_2^n)$, then $\mathfrak{E}_n$ sends 
$u\otimes v$ to $(c_1\otimes c_2 \to V_n((u\otimes v)(c_1\otimes 
c_2)))$, where $c_1 \in C_1$, $c_2 \in C_2$ and $V_n:C_1^n\otimes 
C_2^n \to {(C_1\otimes C_2)}^n$ is the map that shuffles the factors 
of together.  \end{remark}

Now we recall how morphisms of m-coalgebras
were defined in \cite{Smith:1994}:

\begin{definition}
\label{strictmorphdef} Let $C_1$ and $C_2$ be m-coalgebras with sets of
structure maps $\{\ncmap{C_i}_n:C_i \to \homzs{n}(\ncoord{C_i}_n,
C_i^n)\}$, $i = 1, 2$, and all $0\le n < \infty$. A {\it strict
morphism} $\{g,h\}:C_1 \to C_2$ consists of: \begin{enumerate}

\item a chain-map from $g:C_1 \to C_2$;

\item a  morphism of f-resolutions,
$h:\ncoord{C_2} \to \ncoord{C_1}$ such that the diagram

\[
\xymatrix{
C_1\ar[r]^-{\ncmap{C_1}_n}\ar[d]^{g}&\homzs{n}(\ncoord{C_1}_n,C_1^n)
\ar[d]^{\homzs{n}(h,g^n)}\\
C_2\ar[r]_-{\ncmap{C_2}_n}&
\homzs{n}(\ncoord{C_2}_n,C_2^n)}
\]
commutes for all $n$. \end{enumerate} \end{definition}

\begin{definition}
\label{def:scontraction}
A {\it contraction} of chain-complexes 
\[
(f',p,\varphi):C \to D
\] is a
pair of maps $f':C \to D$, $f:D \to C$ and a chain-homotopy $\varphi:C
\to C$ such that:
\begin{enumerate}

\item $f'\circ f = 1_D$

\item $f\circ f' - 1_C =
\partial\varphi$.

\item $\varphi^{2}$=0, $\varphi\circ f=0$, and $f'\circ \varphi=0$ 
\end{enumerate}

The map $f'$ is called the {\it projection} of the
contraction and $f$ is called its {\it injection} --- see
\cite{Gugenheim:1972}.
\end{definition}

\begin{remark}\rem In his thesis 
(\cite{Majewski:1996}), Martin Majewsky called contractions {\it
Eilenberg-Zilber\/} maps.
\end{remark}

\begin{definition} \label{emorphdef} Let $C$ and $D$ be weakly-coherent
m-coalgebras. A contraction 
\[
(f',f,\varphi):C \to D
\]
 with the
injection, $f$, a strict morphism of m-coalgebras, will be called an
{\it elementary equivalence} from $C \to D$. We
will use the notation
\[
\xymatrix{
C\contractright{f}&D
}
\]
to denote an elementary equivalence.
\end{definition} 

\begin{remark}
\rem It is well-known (for instance, see the discussion of Schanuel's 
Lemma in \cite{MacLane:1975}) that any chain-homotopy equivalence of 
two chain-complexes can be decomposed into two iterated contractions.

This implies that contractions are of limited interest when one is
studying {\it chain-complexes.}   This is no longer true when the 
chain-complexes have additional structure --- that of an m-coalgebra, 
for instance.  In this case, the injection of a contraction induces a 
condition on m-structures somewhat similar equivalence of quadratic 
forms.
\end{remark}
 
\begin{definition} \label{mcatdef}  The category of weakly-coherent
m-coalgebras, denoted $\mathfrak{M}$, is defined to be the
localization of $\mathfrak{M}_0$ by the set of strict morphisms whose
associated chain-maps of underlying chain-complexes are injections of
contractions of chain-complexes. \end{definition}

\begin{remarks} \rem The objects of this category are
weakly-coherent m-coalgebras as before, but a morphism from
$A$ to $B$ (say) is a formal composite:
\[
\xymatrix{
{A\ar[r]^-{m_{1}}}
&{\cdots}
&{A_{i}\contractright{s_{i}}}
&*+<20pt>{\cdots}
&{A_{j}\contractleft{s_{j}}}
&*=<-30pt>{\cdots}
&B
}
\]
where the $\{m_{j}\}$ are strict morphisms and the $\{s_{k}\}$ are
elementary equivalences defined in \zref{emorphdef} --- which may go
to the left or right.  We have weakened the definition of morphism
considerably in going from $\mathfrak{M}_{0}$ to $\mathfrak{M}$.
Since projections of contractions are chain-maps, we can still regard
a morphism as having an {\it underlying chain map\/} of chain-complexes.

We will also identify morphisms with the same underlying chain map.

A morphism will be an {\it equivalence\/} if all of its 
constituents are elementary equivalences or their formal inverses.  

\rem The definition is essentially set up so that the maps in the 
Eilenberg-Zilber theorem on page 31 of \cite{Smith:1994} are 
morphisms.  Neither map is a strict morphism, but they both turn out 
to be equivalences.

\rem Morphisms preserve m-structures up to a chain-homotopy.
\end{remarks}

\begin{definition}
\label{def:scellular}
Let $C=\wmco{C}$ be a weakly-coherent m-coalgebra. Then $C$ will be 
called {\it strictly cellular\/} if there exist strict morphisms
of formal coalgebras
\[
g_{k}:\ncoord{C}\to \mathfrak{S}
\]
supporting strict isomorphisms of m-coalgebras
\[
f_{k}:S_{k,n_{k}}=\mathfrak{C}\left(\bigvee_{i=1}^{n_{k}}S^{k-1}
\right)\to C(k-1)
\]
such that 
\[
C(k)=\mathfrak{C}\left(\bigvee_{i=1}^{n_{k}}D^{k}\right)
\bigcup_{f_{k}} C(k-1)
\]
for all $k\ge 0$. Here, $C(k)$ denotes the $k$-skeleton of $C$,
$S_{k,n_{k}}$ is the canonical coherent m-coalgebra of the singular
complex of a wedge of spheres (see 4.2 on page 30 of
\cite{Smith:1994}), and the $D^{k}$ are disks whose boundaries are the
$S^{k-1}$.

We will call a weakly coherent m-coalgebra {\it cellular\/} if it is 
equivalent (in $\mathfrak{M}$) to a strictly cellular m-coalgebra.
\end{definition}

\begin{remarks} \rem If X is a CW complex, 
$\mathfrak{C}(X)=\mathscr{C}(\dot{\Delta}(X))$, where 
$\dot{\Delta}(*)$ is the singular semisimplicial complex functor.

\rem Note that cellularity requires 
the m-structure of an m-coalgebra to be an iterated 
extension of m-structures of spheres.

\rem Clearly, the canonical m-coalgebra of any CW-complex is cellular.  
The converse also turns out to be true --- see 
\zref{cor:realizability}.

It is not hard to find non-cellular m-coalgebras: Consider the
m-coalgebra, $B$, concentrated in dimensions 0 and $3$ (say), where
underlying chain groups are equal to $\integers$.  Equip this with a
trivial coproduct and higher coproducts (subject to the defining
conditions in 3.3 on page 19 of \cite{Smith:1994}). Let $\{e_i\}$
be the generator $\rs{2}$ with boundary $\partial e_i=
(1+{(-1)}^it)e_{i-1}$, where $t\in \integers_2$ is the generator.
We define a map
\[
\Delta:\rs{2}\otimes B \to B\otimes B
\]
where
\begin{enumerate}

\item $B_0=\integers$,

\item $B_3=\integers$, generated by $x$,

\item The ``higher coproducts'' are defined by
\[
\Delta(e_i\otimes x)=\begin{cases}
1\otimes x + x\otimes 1&\text{if $i=0$}\\
0&\text{if $i=2$}\\
x\otimes x&\text{if $i=3$}
\end{cases}
\]
\end{enumerate}
(the last condition is required by 3.3 on page 19 of \cite{Smith:1994}
and implies that the Steenrod operation $Sq^0$ is the identity). Here,
we assume that $t\in \integers_2$ acts trivially on $B$ and multiplies
$B_3\otimes B_3=\integers$ by $-1$.

This is (trivially) coherent -- indeed, it is the m-coalgebra induced
on the {\it homology\/} of the 3-sphere.  It {\it cannot possibly\/}
be cellular because the Hopf invariant of any map from it to a
2-sphere is identically 0.
\end{remarks}

\begin{definition}
\label{def:bigmcat} Define $\hat{\mathfrak{M}}$ be the full 
subcategory of cellular objects of $\mathfrak{M}$.
\end{definition}

We conclude this section with two algebraic results used in the
next section:

\begin{lemma}
\label{lemma:mainlem}
Suppose we have a commutative diagram of weakly-coherent m-coalgebras:
\begin{equation}
\xymatrix{
{A} \ar[d]_{a}\ar@/^3pc/[0,7]^f
&*=<-30pt>{\cdots}
&{U_{i}\contractright{s_{i}}}
&*+<20pt>{\cdots} 
&{U_{j}\contractleft{s_{j}}}
&*=<-30pt>{\cdots} 
&{U_{k}\contractright{s_{k}}}
&{B}\ar[d]^{b}\\
{C}\ar@{=}[0,7]
&*=<-30pt>{}
&
&*+<20pt>{}
&
&*=<-30pt>{}
&
&
{C}
}
\label{dia:firstmap1}
\end{equation}
where the top row is an equivalence from $A$ to
$B$ (whose composite is $f$), and the downward-maps are
strict morphisms.

Then we can expand diagram \zref{dia:firstmap1} to the diagram
\begin{equation}
\xymatrix{
{A} \ar[d]_{a}\ar@/^3pc/[0,7]^f
&*=<-30pt>{\cdots}
&{U_{i}\contractright{s_{i}}\ar[d]_{p_{i}}}
&*+<20pt>{\cdots} 
&{U_{j}\contractleft{s_{j}}\ar[d]_{p_{j}}}
&*=<-30pt>{\cdots} 
&{U_{k}\contractright{s_{k}}\ar[d]_{p_{k}}}
&{B}\ar[d]^{b}
\\
{C}
&*=<-30pt>{\cdots}
&{Z_{i}\contractright{t_{i}}}
&*+<20pt>{\cdots} 
&{Z_{j}\contractleft{t_{j}}}
&*=<-30pt>{\cdots} 
&{Z_{k}\contractright{t_{k}}}
&{C}
}
\label{dia:tworow}
\end{equation}
where
\begin{enumerate}

\item The maps from the first row to the second are all strict morphisms
(see \zref{strictmorphdef}).

\item For all $0 \le i \le k$, the following diagram commutes
\[
\xymatrix{
{U_{i}}\ar[r]^{\varphi_{U_{i}}}\ar[d]_{p_{i}}&{U_{i}}\ar[d]^{p_{i}}\\
{Z_{i}}\ar[r]_{\varphi_{Z_{i}}}&{Z_{i}}
}
\]
where $\varphi_{U_{i}}$ and $\varphi_{Z_{i}}$ are the contracting
homotopies used in the elementary equivalences --- see
\zref{def:scontraction} and \zref{emorphdef}.
\end{enumerate}
\end{lemma}

\begin{proof}
We will actually construct the more complicated diagram:
\begin{equation}
\xymatrix{
{A} \ar[d]_{a}
&*=<-30pt>{\cdots}
&{U_{i}\contractright{s_{i}}\ar[d]_{p_{i}}}
&*+<20pt>{\cdots} 
&{U_{j}\contractleft{s_{j}}\ar[d]_{p_{j}}}
&*=<-30pt>{\cdots} 
&{U_{k}\contractright{s_{k}}\ar[d]_{p_{k}}}
&{B}\ar[d]^{b}
\\
{C}
&*=<-30pt>{\cdots}
&{Z_{i}\contractright{t_{i}}}
&*+<20pt>{\cdots} 
&{Z_{j}\contractleft{t_{j}}}
&*=<-30pt>{\cdots} 
&{Z_{k}\contractright{t_{k}}}
&{C}\\
{C}\ar@{=}[u]
&*=<-30pt>{\cdots}
&{C}\contractup{v_{i}}\ar@{=}[r]
&*+<20pt>{\cdots}
&{C}\contractup{v_{j}}\ar@{=}[l]
&*=<-30pt>{\cdots}
&{C}\contractup{v_{k}}\ar@{=}[r]
&{C}\ar@{=}[u]
}
\label{dia:threerow}
\end{equation}

We construct the lower rows by scanning the upper, from left to right,
and:

\begin{enumerate}
\item
Whenever we encounter a subdiagram of the form
\[
\xymatrix{
U_{i}\contractright{s_{i}}\ar[d]_{p_{i}}&U_{i+1}\ar@{.>}[d]\\
Z_{i}\ar@{.>}[r]&*\txt{?\strut}\\
C\contractup{v_i}&
}
\]
We replace the `?' with the push-out ---
$Z_{i+1}=Z_{i}\oplus U_{i+1} /U_{i}$ (embedded via $(s_{i},p_{i})$) ---
and the appropriate maps. This results in the subdiagram
\[
\xymatrix{
{U_{i}}\contractright{s_{i}}\ar[d]_{p_{i}}
&{U_{i+1}}\ar[d]^{p_{i+1}}\\
Z_{i}\contractright{t_{i}}
&Z_{i+1}\\
C\contractup{v_i}\ar@{=}[r]
&C\contractup{v_{i+1}}
}
\]
where
\begin{enumerate}
\item $p_{i+1}$ and $t_i$
are defined by the canonical property of a push-out
and are strict morphisms of m-coalgebras (see \zref{strictmorphdef}).

\item $t'_i=(1,p_i\circ s'_i):Z_{i}\oplus U_{i+1} /U_{i}\to Z_i$.  
This map is  surjective
since $s'_i$ is, and we have made explicit use of the fact that
$s'_i$ is a left-inverse of $s_i$.

We define a contracting homotopy
$\varphi_{Z_{i+1}}=(0,\varphi_{U_{i+1}}):Z_{i+1}\to Z_{i+1}$, 
where $\varphi_{U_{i+1}}$ is the contracting homotopy of the upper row 
(which exists because it is an elementary equivalence --- 
see \zref{emorphdef}). This makes the lower row an elementary equivalence.

\item $v_{i+1}=(v_i,0):H\to Z_{i+1}=Z_{i}\oplus U_{i+1} /U_{i}$ and
$v'_{i+1}=v'_i\circ t'_i$
\end{enumerate}

\item Whenever we encounter a subdiagram of the form
\[
\xymatrix{
U_{i}\ar[d]_{p_{i}}&U_{i+1}\contractleft{s_{i}}\ar@{.>}[d]\\
Z_{i}&*\txt{?\strut}\ar@{.>}[l]\\
C\contractup{v_i}&
}
\]
we simply pull back $Z_{i}$ to form the diagram
\[
\xymatrix{
U_{i}\ar[d]_{p_{i}}&U_{i+1}\ar[l]_{s_{i}}\ar[d]^{p_{i}\circ s_{i}}\\
Z_{i}&Z_{i}\ar[l]^{1}\\
C\contractup{v_i}\ar@{=}[r]&C\contractup{v_{i+1}}
}
\] 
where $v_{i+1}=v_i$.
\end{enumerate}

This procedure works until we come to the end  (i.e., the 
right end of diagram \zref{dia:threerow}).
\begin{equation}
\xymatrix{
{B}\ar[d]_{\overline{b}}\\
Z_t\\
C\contractup{v_t}
}
\label{dia:rightend}
\end{equation}
where $\overline{b}$ is induced by $b$
--- its target is the embedded copy of $C$.

The commutativity of diagram \zref{dia:firstmap1} implies that we
can splice an extra column onto diagram \zref{dia:rightend} to get
\begin{equation}
\xymatrix{
{B}\ar[d]_{\overline{b}}\ar@{=}[r]
&{B}\ar[d]_{b}\\
Z_t&C\contractleft{v_t}\\
C\contractup{v_t}\ar@{=}[r]&C\ar@{=}[u]
}
\label{dia:newrightend}
\end{equation}
\end{proof}

\begin{corollary}
\label{cor:inductstep}
Suppose we have a commutative diagram of weakly-coherent m-coalgebras:
\begin{equation}
\xymatrix{
{A}\ar[r]^{f}
\ar[d]_{a}
&{B}
\ar[d]^{b}\\
{C}\ar[r]_{=}&{C}
}
\label{dia:firstmap2}
\end{equation}
where the top row is an equivalence, and the downward-maps are
strict morphisms.

Then there exists an equivalence of weakly-coherent m-coalgebras
\[
\hat{f}:A\otimes_{\alpha\circ a}
\cobar C \to
B\otimes_{\alpha\circ b}
\cobar C
\]
where $\cobar (*)$ denotes the cobar construction, $\alpha:C\to \cobar
C$ is the canonical twisting cochain, and the twisted tensor products
are equipped with the canonical weakly-coherent m-structures described
in Proposition 1.19 on page 84 of \cite{Smith:1994}.

In addition, the following diagram commutes:
\begin{equation}
\xymatrix{
{A\otimes_{\alpha\circ a}
\cobar C}\ar[r]^{f\otimes 1}\ar[d]_{1\otimes \epsilon}
&{B\otimes_{\alpha\circ b}
\cobar C}\ar[d]_{1\otimes \epsilon}\\
{A}\ar[r]_{f}
&{B}
}
\label{dia:cobarcommute}
\end{equation}
\end{corollary}

\begin{remark}\rem We will use this and the results of \cite{Smith:1994}
to show that the equivalence $\mathscr{C}(X_1) \to \mathscr{C}(X_2)$
implies the existence of an equivalence between the next stages of
{\it Postnikov towers\/} of $X_1$ and $X_2$. 
\end{remark}

\begin{proof}
This follows by taking diagram \zref{dia:tworow} and putting a third row of
cobar constructions and twisting cochains
\begin{equation}
\xymatrix{
{A} \ar[d]_{a}
&*=<-30pt>{\cdots}
&{U_{i}\contractright{s_{i}}\ar[d]_{p_{i}}}
&*+<20pt>{\cdots} 
&{U_{j}\contractleft{s_{j}}\ar[d]_{p_{j}}}
&*=<-30pt>{\cdots} 
&{U_{k}\contractright{s_{k}}\ar[d]_{p_{k}}}
&{B}\ar[d]^{b}
\\
{C}\ar[d]_{\alpha}
&*=<-30pt>{\cdots}
&{Z_{i}\contractright{t_{i}}}\ar[d]_{\alpha_i}
&*+<20pt>{\cdots} 
&{Z_{j}\contractleft{t_{j}}}\ar[d]_{\alpha_j}
&*=<-30pt>{\cdots} 
&{Z_{k}\contractright{t_{k}}}
\ar[d]_{\alpha_k}
&{C}\ar[d]^{\alpha}\\
{\cobar C}
&*=<-30pt>{\cdots}
&{\cobar Z_{i}\contractright{\cobar(t_{i})}}
&*+<20pt>{\cdots}
&{\cobar Z_{j}\contractleft{\cobar(t_{j})}}
&*=<-30pt>{\cdots}
&{\cobar Z_{k}\contractright{\cobar(t_{k})}}
&{\cobar C}
}
\label{dia:firstcobar}
\end{equation}
where $\alpha_i: Z_{i}\to \cobar  Z_{i}$ are the canonical twisting cochains.

The elementary equivalences on the bottom row are the result of applying 
Proposition 2.32 on page 58 of \cite{Smith:1994}.
\end{proof}

\section{Topological realization of morphisms}
\label{sec:mainresults}

In this section, we will prove the main results involving the
topological realization of m-coalgebras and morphisms.  We begin with
a proof that {\it equivalences\/} topologically realizable
m-coalgebras are topologically realizable.

\begin{theorem}
\label{th:equiv}
Let $X_{1}$ and $X_{2}$ be pointed, simply-connected, locally-finite,
simisimplicial sets, with associated canonical m-coalgebras,
$\mathscr{C}(X_{i})$, $i=1,2$.

In addition, suppose there exists an equivalence of m-coalgebras
\[
f:\mathscr{C}(X_{1}) \to \mathscr{C}(X_{2})
\]
as defined in \cite{Smith:1994} or in \zref{mcatdef} and the
surrounding discussion.

Then there exist refinements (simplicial subdivisions) $X'_i$,
$i=1,2$, of $X_i$, respectively and a simplicial map
\[
\hat{f}:X'_1\to X'_2
\]
such that
\[
f'=\mathscr{C}(\hat{f}):\mathscr{C}(X'_{1}) \to \mathscr{C}(X'_{2})
\]

Consequently, any m-coalgebra equivalence is topologically
realizable up to a chain-homotopy.
\end{theorem}

\begin{remarks}\rem We work in the simplicial category because
the functors $\mathscr{C}(*)$ were originally defined over it. 

It is well-known that the category of locally-finite simplicial
sets coincides with the category of CW complexes. We could also
have worked with the functors $\mathfrak{C}(*)$, computed from
singular complexes.

\rem The refinement is a barycentric subdivision whose degree is
finite within a neighborhood of each vertex of the $X_i$, if they are
finite dimensional.  If the $X_i$ are finite, we can bound this degree
by a finite number.

In any case, however, there are canonical equivalences
\[
\mathscr{C}(X_i)\cong \mathscr{C}(X'_i)
\]
for $i=1,2$.
\end{remarks}

\begin{proof}
The hypothesis implies that the chain-complexes are chain-homotopy 
equivalent, hence that the $X_{i}$, $i=1,2$, have the same homology. 
This implies that the lowest-dimensional nonvanishing homology 
groups --- say $M$ in dimension $k$ --- are isomorphic. We get a 
diagram
\begin{equation}
\xymatrix{
{\mathscr{C}(X_{1})}\ar[r]^{f}
\ar[d]_{\mathscr{C}(c_{1})}
&{\mathscr{C}(X_{2})}
\ar[d]^{\mathscr{C}(c_{2})}\\
{\mathscr{C}(K(M,k))}\ar[r]_{=}&{\mathscr{C}(K(M,k))}
}
\label{dia:firstmap}
\end{equation}

Here, the maps are defined as follows:

\begin{enumerate}

\item The maps $\{\mathscr{C}(c_{i})\}$, $i=1,2$, are induced by
geometric classifying maps;

\item  $f$ is the composite of rightward arrows in the equivalence
between the $\mathscr{C}(X_{i})$, $i=1,2$:
\begin{equation}
\xymatrix{
{\mathscr{C}(X_{1})} 
&*=<-30pt>{\cdots}
&{U_{i}\contractright{s_{i}}}
&*+<20pt>{\cdots} 
&{U_{j}\contractleft{s_{j}}}
&*=<-30pt>{\cdots} 
&{U_{k}\contractright{s_{k}}}
&{\mathscr{C}(X_{2})}
}
\label{dia:mequiv}
\end{equation}
where the $\{U_{\alpha}\}$ are all weakly-coherent m-coalgebras and the
$\{s_{*}\}$  all define elementary equivalences (see 
\zref{emorphdef}). 
\end{enumerate} 

{\bf Claim:} If we forget simplicial structures (i.e., regard the
simplicial sets in \zref{dia:firstmap} as CW-complexes), we may assume
that diagram \zref{dia:firstmap} commutes exactly. to be precise: 

\begin{enumerate}

\item The cellular chain complexes of the $X_i$
are naturally isomorphic to the underlying chain-complexes
of the $\mathscr{C}(X_i)$.

\item We construct the map $c_1$ by finding a topological realization
of the composite $\mathscr{C}(c_2)\circ f$.  That this {\it can\/} be
done follows by elementary obstruction theory and the fact that all
the spaces in question are simply-connected --- see
\cite{Wall:finite2}, for instance. We replace the simplicial map,
$c_1$, by a cellular map, $c'_1$, homotopic to it.The result is a map
of pairs
\[
((X_{1})_{k},(X_{1})_{k-1})\to (X_{2})_{k},(X_{2})_{k-1})
\]
(where $(X_{1})_{k}$ denotes the $k$-skeleton) for all $k\ge 0$, such that
the induced map of cellular chain modules
\[
\pi_{k}((X_{1})_{k},(X_{1})_{k-1})=\mathscr{C}(X_{1})_{k}\to
\pi_{k}((X_{2})_{k},(X_{2})_{k-1})=\mathscr{C}(X_{2})_{k}
\]
exactly coincides with $f$  (regarded only as a map of chain
complexes).

\item Now, we refine the simplicial sets until we can
replace $c'_1$ by a simplicial approximation.  The image of each
simplex of $X_1$ lies in a finite subcomplex of $K(M,1)$ and $X_2$, so
we can simplicially approximate the restriction of $c'_1$ to this
simplex. Consequently, a {\it finite\/} (but, possibly, unbounded)
number of subdivisions of each simplex suffices.
\end{enumerate}

In the following discussion, we will  assume that this
subdivision and simplicial approximation has been carried out
--- and we will suppress the extra notation (i.e., the prime)
for the subdivided complexes and induced maps.

All of the maps in \zref{dia:firstmap} are strict
morphisms of m-coalgebras (see \zref{strictmorphdef}), except for the
map $f$: The vertical maps and the lower horizontal map are strict
because they were induced by {\it geometric\/} maps.

Corollary \zref{cor:inductstep} implies that there exists
an equivalence
\[
\hat{f}:{\mathscr{C}(X_1)}\otimes_{\alpha\circ \mathscr{C}(g_1)}
{\cobar \mathscr{C}(K(M,k))} \to
{\mathscr{C}(X_2)}\otimes_{\alpha\circ \mathscr{C}(g_1)}
{\cobar \mathscr{C}(K(M,k))}
\]
such that the following diagram commutes:
\begin{equation}
\xymatrix{
{\mathscr{C}(X_1)}\otimes_{\alpha\circ \mathscr{C}(g_1)}
{\cobar \mathscr{C}(K(M,k))}\ar[r]^{\hat{f}}\ar[d]_{1\otimes \epsilon}
&{\mathscr{C}(X_2)}\otimes_{\alpha\circ \mathscr{C}(g_1)}
{\cobar \mathscr{C}(K(M,k))}\ar[d]_{1\otimes \epsilon}\\
{\mathscr{C}(X_1)}\ar[r]_{f}
&{\mathscr{C}(X_2)}
}
\label{dia:cobarcommute2}
\end{equation}

Lemma 3.1 of page 93 and Corollary 3.5 on page 96 of \cite{Smith:1994}
imply the existence of equivalences (of weakly-coherent m-coalgebras)
\[
\mathscr{C}(X_i\times_{\hat{\alpha}\circ g_i}
\Omega K(M,k))\to {\mathscr{C}(X_i)}\otimes_{\alpha\circ \mathscr{C}(g_i)}
{\cobar \mathscr{C}(K(M,k))}
\]
for $i=1,2$

We conclude that there is an equivalence
\[
\hat{F}:\mathscr{C}(X_1\times_{\hat{\alpha}\circ g_1}
\Omega K(M,k)) \to
\mathscr{C}(X_2\times_{\hat{\alpha}\circ g_2}
\Omega K(M,k))
\]
where $\Omega (*)$ denotes the loop space functor and
$\hat{\alpha}:K(M,k)\to \Omega K(M,k)$ is the canonical twisting
function (defining a fibration as twisted Cartesian product --- see
\cite{Gugenheim:1959}).

In addition, the commutativity of \zref{dia:cobarcommute2} implies that
\[
f^{*}(\mu_2)=\mu_1\in H^{k+1}(X_1,M)
\]
where $\mu_1$ and $\mu_2$ are the $k$-invariants of the fibrations
$X_1\times_{\hat{\alpha}\circ g_1}
\Omega K(M,k)$ and $X_2\times_{\hat{\alpha}\circ g_2}
\Omega K(M,k)$, respectively.

Since the $X_i\times_{\hat{\alpha}\circ g_i}
\Omega K(M,k)$ are homotopy fibers of the $g_i$ maps for $i=1,2$,
respectively, we conclude that the second stage of the Postnikov
towers of $X_1$ and $X_2$ are equivalent.

A straightforward induction implies
that all {\it finite stages\/} of the Postnikov tower of $X_1$ are
equivalent to corresponding finite stages of the Postnikov tower of
$X_2$. It follows that all finite-dimensional obstructions to 
realizing the underlying chain-map of $f$ by a geometric map
of CW-complexes vanish. 

It is necessary to make one last remark regarding our simplicial
approximations to maps in diagrams like \zref{dia:firstmap} that arise
during inductive steps. Clearly, after any {\it finite\/} number of
inductive steps, we are still dealing with {\it finite\/}
subdivisions of the simplicial sets from the hypothesis. If the
original spaces were finite dimensional, we only need a finite
number of inductive steps.

The conclusion follows.
\end{proof}

Next, we prove a similar result for well-behaved morphisms that aren't
{\it a priori\/} equivalences. We are heading toward a proof that
arbitrary morphisms are topologically realizable.

\begin{prop}
\label{prop:equiv2}
Let $X_{1}$ and $X_{2}$ be pointed, 
simply-connected simisimplicial complexes
complexes, with associated canonical m-coalgebras, 
$\mathscr{C}(X_{i})$, $i=1,2$.

In addition, suppose there exists a strict morphism of
weakly coherent m-coalgebras that induces homology isomorphisms
in all dimensions
\[
f:\mathscr{C}(X_{1}) \to \mathscr{C}(X_{2})
\]
as defined in \cite{Smith:1994}
or in \zref{mcatdef} and the surrounding discussion.

Then there exists a map of CW-complexes (i.e, we forget the
semi-simplicial structure of the spaces and regard them as
CW-complexes --- or pass to suitable simplicial refinements, as
in~\zref{th:equiv}):
\[
\hat{f}:X_1\to X_2
\]
such that
\[
f=\mathscr{C}(\hat{f})
\]

Consequently, $f$ is an equivalence.
\end{prop}

\begin{remarks}\rem This is interesting because strict morphisms don't
generally define m-coalgebra equivalences --- even when they are
homology equivalences. The topological realizability of the
m-coalgebras in question is crucial here.

\rem We could actually have stated that the map $f$ is a composite
$e_1\circ f'\circ e_2$, where $e_1$ and $e_2$ are equivalences
of m-coalgebras and $f'$ is a strict morphism inducing homology
isomorphisms.
\end{remarks}

\begin{proof}
We follow an argument exactly like that used in \zref{prop:equiv2}
above.  In each inductive step we have a morphism of the form
$e_1\circ f_i\circ e_2$, where $e_1$ and $e_2$ are equivalences of
m-coalgebras and $f_i$ is a strict morphism inducing homology
isomorphisms. the only thing we must do differently, here, is to
invoke the Serre Spectral Sequence of a fibration to verify that the
$f_{i+1}$ will be a homology equivalence, given that $f_i$ is.
\end{proof}

\begin{corollary}
\label{cor:homequiv}
Suppose $C_1$ and $C_2$ are weakly coherent m-coalgebras that are
topologically realizable --- i.e., they are equivalent in
$\mathfrak{M}$ (see~\zref{mcatdef}) to $\mathscr{C}(X_i)$,
respectively, for two pointed, simply-connected semi-simplicial
complexes, $X_i$, $i=1,2$.

Then a morphism
\[
f:C_1 \to C_2
\]
is an equivalence if and only if it induces isomorphisms in homology.
\end{corollary}

\begin{theorem}
\label{th:morphism}
Let $X_{1}$ and $X_{2}$ be pointed, simply-connected, locally-finite,
simisimplicial sets complexes, with associated canonical m-coalgebras,
$\mathscr{C}(X_{i})$, $i=1,2$.

In addition, suppose there exists a  morphism of m-coalgebras
(see~\zref{strictmorphdef}):
\[
f:\mathscr{C}(X_{1}) \to \mathscr{C}(X_{2})
\]

Then there exists a map of CW-complexes (i.e, we forget the
semi-simplicial structure of the spaces and regard them as
CW-complexes --- or form simplicial refinements, as
in~\zref{th:equiv}):
\[
\hat{f}:X_1\to X_2
\]
such that
\[
f=\mathscr{C}(\hat{f})
\]

Consequently, any morphism of m-coalgebras is topologically realizable
up to a chain-homotopy.
\end{theorem}

\begin{proof}
We prove this result by an inductive argument somewhat different from
that used in theorem~\zref{th:equiv}.

We  build a sequence of fibrations
\[
\xymatrix{
F_i\ar[d]_{p_i}\\
X_2
}
\]
over $X_2$ in such a way that
\begin{enumerate}

\item the morphism $f:\mathscr{C}(X_1)\to \mathscr{C}(X_2)$ lifts to
$\mathscr{C}(F_i)$ --- i.e., we have commutative diagrams
\[
\xymatrix{
&{\mathscr{C}(F_i)}\ar[d]^{\mathscr{C}(p_i)}\\
{\mathscr{C}(X_1)}\ar[r]_{f}\ar[ur]^{f_i}&{\mathscr{C}(X_2)}
}
\]

For all $i>0$, $F_i$ will be a fibration over $F_{i-1}$ with fiber
a suitable Eilenberg-MacLane space.

\item The map $f_i$ is $i$-connected in homology.
\end{enumerate}

If the morphism $f$ were geometric, we would be building
its {\it Postnikov tower.}

Assuming that this inductive procedure can be carried out, we note
that it forms a convergent sequence of fibrations (see
\cite{Spanier:1966}, chapter 8, \S~3).  This implies that we may pass
to the inverse limit and get a commutative diagram
\[
\xymatrix{
&{\mathscr{C}(F_{\infty})}\ar[d]^{\mathscr{C}(p_{\infty})}\\
{\mathscr{C}(X_1)}\ar[r]_{f}\ar[ur]^{f_{\infty}}&{\mathscr{C}(X_2)}
}
\]
where $\tilde{f}_{\infty}$ is a morphism of weakly-coherent m-coalgebras
that is a {\it homology equivalence.} Now \zref{cor:homequiv}
implies that $f_{\infty}$ is an equivalence of m-coalgebras, and
\zref{th:equiv} implies that it is topologically realizable.

It follows that we get a (geometric) map
\[
\bar{f}_{\infty}:X_1 \to F_{\infty}
\]
and the composite of this with the projection 
$p_{\infty}:F_{\infty}\to X_2$ is a topological realization
of the original map $f:\mathscr{C}(X_1)\to \mathscr{C}(X_2)$.

It only remains to verify the inductive step:

Suppose we are in the $k^{\text{th}}$
iteration of this inductive procedure. Then
the mapping cone, $\mathscr{A}(f)$ is
acyclic below dimension $k$. Suppose  that
$H_k(\mathscr{A}(f_k))=M$.  Then we get a long exact sequence in
cohomology:
\begin{multline}
{\dots}\to H^{k+1}(X_1,M)\to H^k(\mathscr{A}(f_k),M)=\homz(M,M)\\
\to H^k(F_k,M)\to H^k(X_1,M)\to 0
\label{fib:exact}
\end{multline}

Let $\mu\in H^k(F_k,M)$ be the image of $1_M\in
H^k(\mathscr{A}(f_k),M)=\homz(M,M)$ and consider the map
\[
h_{\mu}:X_2 \to K(M,k)
\]
classified by $\mu$.  We pull back the contractible fibration
\[
K(M,k)\times_{\bar{\alpha}} \Omega K(M,k)
\]
over $h_{\mu}$ to get a fibration
\[
F_{k+1}=F_k\times_{\bar{\alpha}\circ h_{\mu}}\Omega K(M,k)
\]
where, as before, $\Omega(*)$ represents the loop space.

{\bf Claim:} The morphism $f_k$ lifts to a morphism
$f_{k+1}:\mathscr{C}(X_1)\to \mathscr{C}(F_{k+1})$ in such
a way that the following diagram commutes:
\[
\xymatrix{
{\mathscr{C}(X_1)}\ar[r]^-{f_{k+1}}\ar[rd]_{f_k}
&{\mathscr{C}(F_{k+1})}\ar[d]^{p}\\
&{\mathscr{C}(F_k)}
}
\]
where $p'_{k+1}:F_{k+1}\to F_k$ is that fibration's projection map.

{\bf Proof of Claim:} We begin by using Lemma 3.1 of page 93 and
Corollary 3.5 on page 96 of \cite{Smith:1994} to conclude the
existence of a commutative diagram:
\begin{equation}
\xymatrix{
{\mathscr{C}(F_k\times_{\bar{\alpha}\circ h_{\mu}}\Omega K(M,k))}
\ar[r]^{e}\ar[dr]_{p'_{k+1}}
&{\mathscr{C}(F_k)
\otimes_{\bar{\alpha}\circ h_{\mu}}\cobar \mathscr{C}(K(M,k))}
\ar[d]^{1\otimes \epsilon}\\
&{\mathscr{C}(F_k)}
}
\label{dia:timesequiv}
\end{equation}
where $e$ is an m-coalgebra equivalence.

If we pull back this twisted tensor product over the map $f_k$, we get
a trivial twisted tensor product (i.e., an untwisted tensor product),
because the image of $f^*(\mu)=0\in H^{k}(X_1,M)$, by the exactness of
\zref{fib:exact}. Theorem 1.20 on page 85 of \cite{Smith:1994}
implies the existence of a morphism
\begin{multline}
{\mathscr{C}(X_1)}\to {\mathscr{C}(X_1)}\otimes 1\subset 
{\mathscr{C}(X_1)}\otimes \cobar {\mathscr{C}(K(M,k))}\\ \to
{\mathscr{C}(F_k)}\otimes_{\bar{\alpha}\circ h_{\mu}}
\cobar{\mathscr{C}(K(M,k))}
\end{multline}
The composition of this map with $e$ in \zref{dia:timesequiv} is the
required map
\[
{\mathscr{C}(X_1)}\to {\mathscr{C}(F_{k+1})}
\]

To see that $H_{k}(\mathscr{A}(f_{k+1}))=0$, note that:
\begin{enumerate}

\item $\mu\in H^k(F_k,M)= H^k({\mathscr{C}(F_k)},M)$
is the pullback of the class in 
$H^{k}(\mathscr{A}(f_k),M)$ inducing a homology isomorphism
\[
\mu:H_{k}(\mathscr{A}(f_k))\to H_{k}(K(M,k))
\]
(by abuse of notation, we identify $\mu$ with a cochain) or
\[
\mu:H_{k}(\mathscr{A}(f_k))\to H_{k}({\mathscr{C}(K(M,k))})
\]

\item in the stable range,  ${\mathscr{C}(F_k)}
\otimes_{\bar{\alpha}\circ h_{\mu}}
\cobar{\mathscr{C}(K(M,k))}$ is nothing but the {\it algebraic mapping 
cone\/} of the chain-map, $\mu$, above.  But the algebraic mapping cone 
of $\mu$ clearly has vanishing homology in dimension $k$ since $\mu$ 
induces homology isomorphisms.
\end{enumerate}
\end{proof}

\begin{corollary}
\label{cor:realizability} A weakly-coherent m-coalgebra
is topologically realizable if and only if it is cellular
(see \zref{def:scellular}).
\end{corollary}

\begin{proof}
Clearly, topologically realizable m-coalgebras are cellular.

Theorem \zref{th:morphism} implies the converse, because all of the 
attaching morphisms in \zref{def:scellular} are topologically realizable. 
\end{proof}

\begin{corollary}
\label{cor:mainresult} The functor 
\[
\mathscr{C}(*):\underline{\mathrm{Homotop}}_0\to \hat{\mathfrak{M}}
\]
(see \zref{def:bigmcat} for the definition of $\hat{\mathfrak{M}}$)
defines an equivalence of categories, where
$\underline{\mathrm{Homotop}}_0$ is the category of pointed,
simply-connected CW-complexes and continuous maps, in which homotopy
equivalences have been inverted (i.e., it is the category of fractions
by homotopy equivalences).
\end{corollary}

\bibliographystyle{amsplain} \bibliography{mrabbrev,refs}

\vspace{3mm}{\small\noindent{\sc Department of
Mathematics and Computer Science\\
Drexel University\\
Philadelphia,
PA 19104}
\par
\vspace{3mm}
\noindent Email: \verb,jsmith@mcs.drexel.edu,\\
Home page: \verb,http://www.mcs.drexel.edu/~jsmith,}
\end{document}